# Global well-posedness and exponential decay for a fluid-structure model with small data

Igor Kukavica and Wojciech S. Ożański

Tuesday 30th August, 2022

**Abstract**

We address the system of partial differential equations modeling the motion of an elastic body interacting with an incompressible fluid via a free interface. The fluid is modeled by the incompressible Navier-Stokes equations while the structure is represented by a damped wave equation $w_{tt} - \Delta w + \alpha w_t = 0$, where $\alpha > 0$. We prove the global existence and exponential decay of strong solutions for small initial data in a suitable Sobolev space. We show that the elastic velocity $w_t$ and the acceleration $w_{tt}$ can be controlled by the $H^2$ elliptic estimates and by the dissipation of the fluid via the free interface. We also find that, even though the vanishing of the final displacement $w$ appears invisible in the energy method, it can be deduced from the preservation of total volume. Our approach allows for any superlinear perturbation of the wave equation.

## 1   Introduction

We are concerned with a model of an elastic body interacting with a viscous incompressible fluid. The elastic body occupies a time-dependent domain $\Omega_e(t)$, and is described by a displacement function $w(x,t) := \eta(x,t) - x$, where $\eta(x,t)$ denotes the Lagrangian mapping of a particle located at $x \in \Omega_e(0) =: \Omega_e$. The displacement $w$ is governed by the linear damped equation

$$w_{tt} - \Delta w + \alpha w_t = 0 \tag{1.1}$$

in $\Omega_e \times (0, \infty)$. We note that all the results presented below extend to the case of an additional superquadratic term $f(w)$, i.e., a smooth nonlinearity $f$ such that $f(x), f'(x), f''(x) = o(1)$ as $x \to 0$. For the boundary of $\Omega_e(t)$, we suppose that $\partial \Omega_e(t) = \Gamma_e(t) \cup \Gamma_c(t)$. On the elastic boundary $\Gamma_e(t) = \Gamma_e(0)$ we assume the zero displacement condition

$$w(x,t) = 0, \qquad x \in \Gamma_e, \quad t > 0, \tag{1.2}$$

and $\Gamma_c(t)$ is assumed to be a common boundary shared with a fluid described by the incompressible Navier-Stokes equations

$$u_t - \Delta u + (u \cdot \nabla) u + \nabla p = 0,$$

$$\operatorname{div} u = 0 \qquad \text{in } \{(y,t)\colon t > 0,\, y \in \Omega_f(t)\}$$

([Te1, Te2]), where $\Omega_{\mathrm{f}}(t)$ denotes the time-dependent fluid domain, which, except for the common free interface $\Gamma_{\mathrm{c}}(t)$, has a fixed boundary $\Gamma_{\mathrm{f}} \coloneqq \{x_3 = 2\}$. We consider a simplified form of $\Omega_{\mathrm{e}} = \Omega_{\mathrm{e}}(0)$ and $\Omega_{\mathrm{f}} = \Omega_{\mathrm{f}}(0)$ by assuming periodicity in the variables $x_1, x_2$, and that the fluid is located above the elastic solid. Namely, we take

$$\Omega_{\mathrm{e}} \coloneqq \left\{ x = (x_1, x_2, x_3) \colon (x_1, x_2) \in \mathbb{T}^2,\, 0 \leq x_3 \leq 1 \right\}, \qquad (1.3)$$
$$\Omega_{\mathrm{f}} \coloneqq \left\{ x = (x_1, x_2, x_3) \colon (x_1, x_2) \in \mathbb{T}^2,\, 1 \leq x_3 \leq 2 \right\};$$

see Figure 1 for a sketch.

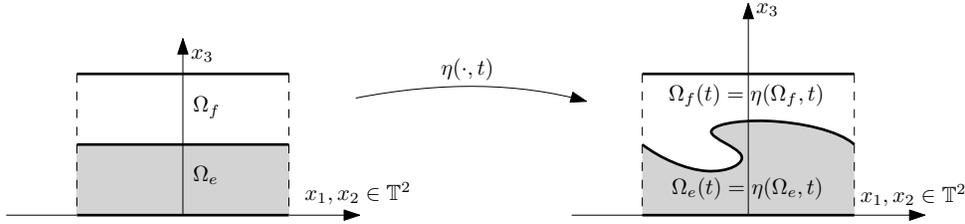

Figure 1: The sketch of the fluid-structure interaction model.

In order to state our main result, we first introduce the Lagrangian setting and specify the boundary conditions on $\Gamma_{\mathrm{c}}$ and $\Gamma_{\mathrm{f}}$.

In the fluid domain $\Omega_{\mathrm{f}}$, the displacement $\eta(\cdot, t)$ is defined as the solution of the system

$$\begin{aligned} \eta_t(x, t) &= v(x, t), \\ \eta(x, 0) &= x, \end{aligned} \qquad (1.4)$$

for $x \in \Omega_{\mathrm{f}}$ and $t > 0$, where

$$v(x, t) \coloneqq u(\eta(x, t), t). \qquad (1.5)$$

The incompressible Navier-Stokes equations (1.3) in the Lagrangian coordinates read

$$\partial_t v_i - \partial_j(a_{jl} a_{kl} \partial_k v_i) + \partial_k(a_{ki} q) = 0 \ \text{ in } \Omega_{\mathrm{f}} \times (0, T), \qquad i = 1, 2, 3, \qquad (1.6)$$
$$a_{ki} \partial_k v_i = 0 \ \text{ in } \Omega_{\mathrm{f}} \times (0, T), \qquad (1.7)$$

where $a \coloneqq (\nabla \eta)^{-1}$ and $q(x, t) \coloneqq p(\eta(x, t), t)$ denotes the Lagrangian pressure. Except for the zero displacement condition (1.2) on the elastic boundary, we also assume continuity of the velocities

$$w_t = v \ \text{ on } \Gamma_{\mathrm{c}} \times (0, T) \qquad (1.8)$$

and continuity of the stresses

$$\partial_j w_i N_j = a_{jl} a_{kl} \partial_k v_i N_j - a_{ji} q N_j \ \text{ on } \Gamma_{\mathrm{c}} \times (0, T), \qquad i = 1, 2, 3, \qquad (1.9)$$

where we use the summation convention on repeated indices and $N = e_3$ denotes the outward unit normal vector with respect to $\Omega_{\mathrm{e}}$. We note that the boundary conditions (1.8)–(1.9) are specified on the common boundary $\Gamma_c$ in the Lagrangian variables, and thus $\Gamma_c$ should be thought of fixed and equal $\Gamma_c(0)$. On the outside fluid boundary $\Gamma_f$, we assume the non-slip boundary condition

$$v = 0 \text{ on } \Gamma_f \times (0, T). \tag{1.10}$$

Note that, for the sake of simplicity, we use the strain tensor $\nabla v$ in (1.9) instead of the symmetric gradient matrix $\nabla v + \nabla v^T$. We thus seek a solution $(v, w, q, a, \eta)$ to the damped fluid-structure system (1.1)–(1.10).

The existence of local-in-time solutions to the system in question was first studied, without damping, by Coutand and Shkoller [CS1, CS2]. Other local-in-time well-posedness results were established for more general initial data in [B, KT1, KT2, RV, BGT]. Solutions for an analogous model involving the compressible Navier-Stokes equations were also considered in [BG1, BG2, KT3]; see also [AL, ALT, AT1, B, BGLT1, BGLT2, BZ1, DEGL, DGHL, GGCC, IKLT1, MC1, MC2] for other results on local strong and global weak solutions for related models.

In [IKLT2] the authors considered the above model without the homogeneous Dirichlet boundary condition (1.2), but with additional *interior damping* of the form of an additional term $\beta w$ on the left-hand side of (1.1), and with an additional stabilizing term $-\gamma \partial_3 w$ on the right-hand side of (1.8), the so-called *transmission boundary conditions*. It was shown that such stabilized system is globally well-posed for small initial data. A subsequent work [IKLT3] considered the global well-posedness for the model without the boundary stabilization (i.e. with $\gamma = 0$), but still included the interior one. The purpose of this work is to address the well-posedness of the model without any stabilization terms.

One of the most remarkable features of the model (1.1)–(1.10) is relatively little control of the lowest order terms. In fact, one could expect that, during the evolution in time of the system, the elastic structure could move to a non-trivial finite state, a phenomenon that appears very difficult to capture. In fact, the role of the extra stabilization terms considered by [IKLT2, IKLT3] was to ensure that this does not happen.

It turns out that, despite the Dirichlet boundary condition on $\Gamma_e$ (which allows us to use the Poincaré inequality to control $\|w\|$ in terms of $\|\nabla w\|$), the energy methods, in fact, do not give sufficient control of the lowest order terms that would imply global well-posedness and exponential decay. To be more precise, it appears that the quantities

$$\int_{\Gamma_c} w \quad \text{and} \quad \int_{\Gamma_c} q \tag{1.11}$$

cannot be controlled by the natural notion of the energy of the system, even if it is assumed to decay to 0.

In order to comment on this phenomenon, we note that we will be using the notation

$$\|\cdot\|_k := \|\cdot\|_{H^k(\Omega)}, \qquad \|\cdot\| := \|\cdot\|_{L^2(\Omega)},$$

where $\Omega = \Omega_f$ or $\Omega = \Omega_e$, depending whether the function is defined in the fluid or elastic domain. We also denote by $\partial'$ the gradient with respect to the variables $x_1$ and $x_2$ only, and by $\partial''$ the matrix of all second order derivatives with respect to $x_1$ and $x_2$.

The dynamics of the system (1.1)–(1.10) can be studied using various notions of energies, and each such notion corresponds to a differential operator $S$ that can be applied to the system. Such an operator must be well-adapted to the system in the sense that it could be applied to both the Navier-Stokes equation (1.6) and the structure equation (1.1) and that it could be transferred through the common boundary. Namely, taking $S$ of (1.6) and (1.1) and testing against appropriate test functions, one must be able to match the resulting boundary terms at $\Gamma_c$ and using the boundary conditions (1.8)–(1.9), see (3.1) and (3.9) for example, which we use to obtain an energy inequality (3.10) for each $S$. It turns out that in such an approach $S$ can involve only derivatives $\partial'$ in the horizontal variables $x' = (x_1, x_2)$ or time derivatives $\partial_t$, namely

$$S \in \{\mathrm{id}, \partial', \partial_t, \partial'\partial_t, \partial'', \partial_{tt}\}, \tag{1.12}$$

which is one of the main structural properties of the problem. In particular, one cannot consider $S = \partial_{x_3}$, as it is not clear how to carry it across the common boundary $\Gamma_c$.

Thus, let us suppose that we can estimate all the terms of the form

$$\|Sv\|, \|Sw\|, \|\nabla Sw\|, \|Sw_t\| \tag{1.13}$$

for $S \in \{\mathrm{id}, \partial', \partial_t, \partial'\partial_t, \partial'', \partial_{tt}\}$, which naturally arise in this approach (see (3.14)). Note that the only way to control $\int_{\Gamma_c} q$ is to employ (1.9) and try to estimate

$$\int_{\Gamma_c} a_{3l} a_{kl} \partial_k v \quad \text{and} \quad \int_{\Gamma_c} \partial_3 w. \tag{1.14}$$

The former can be controlled using a trace estimate by $\|v\|_2$ and then a Stokes-type estimate (see (2.17)) gives a bound in terms of $\|v_t\| + \|v\|_{H^{3/2}(\Gamma_c)}$. The term $\|v_t\|$ is now the first term in (1.13) for $S = \partial_t$, while $\|v\|_{H^{3/2}}$ can be estimated by the third term in (1.13) (for $S \in \{\partial_t, \partial'\partial_t\}$) using the boundary condition (1.8) and a trace estimate,

$$\|w_t\|_{H^{3/2}(\Gamma_c)} \lesssim \|w_t\|_{H^{1/2}(\Gamma_c)} + \|\partial' w_t\|_{H^{1/2}(\Gamma_c)} \lesssim \|w_t\|_1 + \|\partial' w_t\|_1.$$

As for the second term in (1.14) one can use trace estimates to bound it by $\|w\|_2$, which in turn can be bounded by the elliptic estimate with the Dirichlet boundary condition (see (2.12)),

$$\|w\|_2 \lesssim \|\Delta w\| + \|w\|_{H^{3/2}(\Gamma_c)} \lesssim \|w_t\| + \|w_{tt}\| + \|w\|_{H^{3/2}(\Gamma_c)}.$$

The first two terms correspond to the second term in (1.13) for $S \in \{\partial_t, \partial_{tt}\}$. Moreover, the highest order part of the last term, $\|w\|_{H^{3/2}(\Gamma_c)}$ can also be controlled, since a trace estimate and the Poincaré inequality give $\|\partial' w\|_{H^{1/2}(\Gamma_c)} \lesssim \|\partial' w\| + \|\nabla \partial' w\|$, which gives the second and third terms in (1.13) (for $S = \partial'$). However, the lowest order part of $\|w\|_{H^{3/2}(\Gamma_c)}$ requires control of $\int_{\Gamma_c} w$, which is not clear. In an analogous way, one may start with trying to estimate $\int_{\Gamma_c} w$ and end up with the need to estimate $\int_{\Gamma_c} q$. In this sense the control of the two terms (1.14) is dual.

On the other hand, one might think of using a trace estimate and the Poincaré inequality to estimate $\int_{\Gamma_c} w$ by the second and third terms in (1.13) for $S = \text{id}$. However, to this end one needs to use the energy estimate of the system for $S = \text{id}$, which in turn requires control of $\int_{\Gamma_c} q$, and so one obtains a circular argument. In order to see it, we explain below (in Section 3.1) that such energy estimate on time interval $(\tau, t)$ will include a term of the form

$$\int_\tau^t \int_{\Omega_f} Sq \operatorname{div} \widetilde{Sw}(\tau), \tag{1.15}$$

for $S = \text{id}$, where $\widetilde{Sw}$ denotes the Sobolev extension of $Sw$ to the fluid domain $\Omega_f$ (see (3.4)). It turns out that, in order to handle this term one does need to control the lowest order parts of both $q$ and $w$, including $\int_{\Gamma_c} w$ and $\int_{\Gamma_c} q$, which then becomes a circular argument.

We note in passing that the analysis of the term (1.15) of the energy inequalities, for various $S$, is one of the most challenging issues of the system, and we refer the reader to (3.18)–(3.21) for further insights. In either case, at this point it is not clear how to control the terms (1.11), and so also not clear how to proceed.

Remarkably, one can observe that the equations (1.1)–(1.10) of the system suggest that $\int_{\Gamma_c} w$ does not decay to zero. Indeed, using the Fundamental Theorem of Calculus in $x_3$ and the homogeneous Dirichlet boundary condition (1.2), we see that

$$\int_{\Gamma_c} w = \int_{\Omega_e} \partial_3 w =: c. \tag{1.16}$$

Moreover, taking $\partial_3$ of (1.1) and integrating over $\Omega_e$ we obtain an ODE for $c(t)$,

$$c'' + \alpha c' = \int_{\Omega_e} \Delta w,$$

with the initial condition $c(0) = 0$. This gives that

$$c(t) = \int_0^t \int_0^s \int_{\Omega_e} \partial_3 \Delta w(x, h) \mathrm{d}x \, \mathrm{e}^{\alpha(h-s)} \mathrm{d}h \, \mathrm{d}s + c'(0) \left(1 - \mathrm{e}^{-\alpha t}\right),$$

which may not decay to zero even if $c'(0) = 0$ and if $\|D^3 w\|$ decays to zero. However, we show below (in (1.19)) that it does decay to zero merely thanks to the preservation of total volume.

The main result of the paper (Theorem 1 below) uses this observation and considers the *double-normalized wave displacement*, namely

$$\overline{w} := w - x_3 \fint_{\Omega_e} \partial_3 w, \tag{1.17}$$

in order to obtain a global-in-time well-posedness and exponential decay.

**Theorem 1** (Main result)**.** *Let $(v, w, q, \eta, a)$ be a smooth solution to (1.1)–(1.10) on some time interval $[0, T)$, and set*

$$Y(t)^2 \coloneqq \|v\|_3^2 + \|v_t\|_1^2 + \|v_{tt}\|^2 + \|\overline{w}\|_3^2 + \|w_t\|_2^2 + \|w_t\|_2^2 + \|w_{tt}\|_1^2 + \|w_{ttt}\|^2.$$

*Then there exists $C \geq 1$ and $\varepsilon > 0$, independent of $T$, such that if $Y(0) \leq \varepsilon$ then*

$$Y(t) \leq C\varepsilon \mathrm{e}^{-t/C},$$

*for $t \in [0, T)$.*

As in [IKLT3], the a priori estimate of Theorem 1 can be used to show the global-in-time existence and exponential decay of solutions to (1.1)–(1.10) for small initial data.

**Theorem 2** (Global-in-time well-posedness for small data)**.** *There exists $C \geq 1$ such that for every sufficiently small $\varepsilon > 0$ and*

$$\|v_0\|_4, \|\partial_t w_0\|_2, \|\partial_t q_0\|_1, \|q_0\|_3 \leq \varepsilon,$$

*with appropriate compatibility conditions (i.e. [IKLT3, (2.9)–(2.11), (2.15)]), then there exists a unique solution $(v, w, q, \eta, a)$ of (1.1)–(1.10) such that*

$$\|v\|_3 + \|\overline{w}\|_3 + \|w_t\|_2 + \|w_{tt}\|_1 + \|w_{ttt}\| + \|v_{tt}\| \leq C\varepsilon \mathrm{e}^{-t/C},$$

*for $t \geq 0$.*

We note that at first sight the correction (1.17) seems unlikely to solve the problem discussed above, since the PDE for $\overline{w}$ is

$$\overline{w}_{tt} + \alpha \overline{w}_t = \Delta \overline{w} - x_3 d,$$

where

$$d(t) \coloneqq \fint_{\Omega_\mathrm{e}} \Delta \partial_3 w(t).$$

This is different than the equation (1.1) for $w$, and so one could expect further obstacles to arise. However, $\overline{w}$ enjoys some cancellation properties. First,

$$\int_{\Gamma_\mathrm{c}} \overline{w} = \int_{\Gamma_\mathrm{c}} w - |\Gamma_\mathrm{c}| \fint_{\Omega_\mathrm{e}} \partial_3 w = \int_{\Gamma_\mathrm{c}} w - \int_{\Gamma_\mathrm{c}} w = 0, \tag{1.18}$$

where we recalled (1.16) and used the fact that $|\Gamma_\mathrm{c}| = |\Omega_\mathrm{e}|$. Thus, noting that $\int_{\Gamma_\mathrm{c}} \partial' \overline{w} = 0$ as well, we see that $\partial' \overline{w} = \partial' w$ (recall (1.17)), and so the elliptic regularity (see (2.12) below) for the equation $\Delta \overline{w} = \Delta w$ gives

$$\|\overline{w}\|_2 \lesssim \|\Delta \overline{w}\| + \|\overline{w}\|_{H^{3/2}(\Gamma_\mathrm{c})} \lesssim \|\Delta \overline{w}\| + \|\partial' \overline{w}\|_{H^{1/2}(\Gamma_\mathrm{c})}$$
$$= \|\Delta w\| + \|\partial' w\|_{H^{1/2}(\Gamma_\mathrm{c})} \lesssim \|w_{tt}\| + \|w_t\| + \|\partial' w\|_1,$$

which can now we bounded by the terms (1.13). Furthermore, one expects the viscosity of the fluid to stabilize the dynamics of the elastic structure. This is indeed the case, and such stabilization can be observed in some of our estimates in Section 3 below. In fact, thanks to the introduction of the correction $\overline{w}$, we do not need to consider the case $S = \mathrm{id}$ (recall (1.12)) in our energy estimates which would have caused difficulties due to the problematic term (1.15) discussed above.

As explained above, the well-posedness result of Theorem 2 is not expected to yield an exponential decay of the quantities (1.11). However, it can be shown that under the conditions of Theorem 2, they remain of order $\varepsilon$ for all times, see (2.23) and (3.17). Moreover, one can see that the final state is zero by the preservation of total volume. Indeed, we first note that Theorem 2 implies that $w(t) \to x_3 W$ for some $W \in \mathbb{R}^3$ with $|W| \lesssim \varepsilon$. Thus

$$\lim_{t \to \infty} |\Omega_\mathrm{e}(t)| = \int_{\Omega_\mathrm{e}} |\det(\nabla(x + x_3 W))| \mathrm{d}x = \pi^2 |1 + W_3|.$$

Hence, since $|\Omega_\mathrm{f}(t)| = |\Omega_\mathrm{f}(0)| = \pi^2$ for all times, due to the incompressibility of the fluid, and $|\Omega_\mathrm{f}(t)| + |\Omega_\mathrm{e}(t)| = 2\pi^2$ as $t \to \infty$, which follows from $w|_{\Gamma_\mathrm{e}} = 0$ and $v|_{\Gamma_\mathrm{f}} = 0$, we obtain

$$\pi^2 = \lim_{t \to \infty} (2\pi^2 - |\Omega_\mathrm{f}(t)|) = \lim_{t \to \infty} |\Omega_\mathrm{e}(t)| = \pi^2 |1 + W_3|, \tag{1.19}$$

which implies that $W_3 = 0$. Therefore, taking $t \to \infty$ in the third equation of (1.9) (i.e., for $i = 3$) gives that $\lim_{t \to \infty} \int_{\Gamma_\mathrm{c}} q = 0$, which in turn implies that $W_1 = W_2 = 0$, by taking $t \to \infty$ in the first and second equations of (1.9). This shows that the correction in (1.17) decays to zero as $t \to \infty$, but note that it is not clear whether it vanishes at any other $t > 0$.

Except for the strategy described above, the proof of Theorem 1 follows a similar scheme as in [IKLT2, IKLT3], but involves a number of improvements and new ideas, such as immense simplifications of all the estimates, an utilization of a new test function $\phi$ (see (3.3)), the double-normalized wave displacement (1.17), and also avoids the need for the lowest order energy estimate. Namely, for each $S \in \{\partial', \partial_t, \partial'\partial_t, \partial'', \partial_{tt}\}$ we derive an energy inequality (see (3.10)) that expresses an energy coupling between the fluid and the elastic solid that also depends on a coupling between two energy levels (see (3.1) and (3.9)). The coupling is quantified using a parameter $\lambda > 0$, which then needs to be taken sufficiently small.

Using such energy inequalities, we then select some quantities appearing in the inequalities to identify a total energy $X(t)$ of the system (as in (1.13), see (3.14) below). We then estimate the behavior in time of all ingredients of the total energy (see Steps 1–3 in Section 3.2) in a way that enables us to fix $\lambda$ sufficiently small (see (3.24)) to employ an ODE-type lemma (Lemma 4) to deduce exponential decay of the total energy for small initial data.

The paper is structured as follows. In Section 2 we introduce notation and discuss some properties of the Lagrangian map $\eta$ as well as the Stokes estimates concerning (1.6). We then prove our main result, Theorem 1, in Section 3. The proof is based on a number of energy estimates, given in Section 3.1, which lead to the definition (3.14) of the total energy $X(t)$ in Section 3.2, where we also combine the energy estimates into an a priori bound (3.23) on $X$. The exponential decay of $X$ is then established using an ODE-type result, Lemma 4, which is proven in Section 3.3.

## 2 Preliminary results

Throughout the paper we are not concerned with the dependence of our estimates on $\alpha > 0$. In particular all constants, which are denoted by $C$, may depend on $\alpha$. The values of constants $C$ may change from line to line. We apply the summation convention over repeated indices.

### 2.1 Estimates of the particle map

We note that, for each $t > 0$, we have

$$\eta = x + \int_0^t v, \tag{2.1}$$

due to (1.4), and so in particular

$$\|I - \nabla \eta(t)\|_2 + \|D^2 \eta(t)\|_1 \leq \|v\|_{L^1((0,t);H^3)}, \tag{2.2}$$

where $I$ denotes the $3\times 3$ identity matrix. for every $t > 0$. Moreover, due to the incompressibility condition in (1.7), we have that $\det \nabla \eta = 1$ for all times, which shows that $a = (\nabla \eta)^{-1}$ is the corresponding cofactor matrix, that is,

$$a_{ij} = \frac{1}{2}\epsilon_{imn}\epsilon_{jkl}\partial_m \eta_k \partial_n \eta_l, \tag{2.3}$$

where $\epsilon_{ijk}$ denotes the permutation symbol. By (2.2) and (2.3), we have

$$\begin{aligned}\|a_t\|_2 &\lesssim \|\nabla \eta\|_2 \|v\|_3 \lesssim \|v\|_{L^\infty((0,t);H^3)}(1 + \|v\|_{L^1((0,t);H^3)}),\\ \|a_{tt}\|_1 &\lesssim \|\nabla \eta\|_2 \|v_t\|_2 + \|v\|_3^2\end{aligned} \tag{2.4}$$

and

$$\|I - a\|_2 \leq \int_0^t \|a_t\|_2 \lesssim \int_0^t \|\nabla \eta\|_2 \|v\|_3 \lesssim \|v\|_{L^1((0,t);H^3)}\left(1 + \|v\|_{L^1((0,t);H^3)}\right), \tag{2.5}$$

for every $t > 0$. Moreover,

$$\|I - aa^T\|_2 \lesssim \|I - a\|_2 + \|a\|_2 \|I - a^T\|_2 \lesssim \|I - a\|_2 \left(1 + \|I - a\|_2\right), \tag{2.6}$$

and so altogether we obtain that for every sufficiently small $\gamma > 0$ we have

$$\|a_t\|_2 + \|I - a\|_{L^\infty} + \|I - aa^T\|_{L^\infty} + \|I - a\|_2 + \|I - aa^T\|_2 + \|I - \nabla \eta\|_2 \lesssim \gamma, \tag{2.7}$$

for all $t \geq 0$ such that $(\|v\|_{L^1((0,t);H^3)} + \|v\|_{L^\infty((0,t);H^3)})(1 + \|v\|_{L^1((0,t);H^3)}) \leq \gamma$, where we have again employed (2.2)–(2.3). We also recall the Piola identity

$$\partial_i a_{ij} = 0, \tag{2.8}$$

for all $i = 1, 2, 3$, which can be verified directly by (2.3).

## 2.2  Stokes-type and elliptic estimates

In the first statement, we recall the classical Sobolev regularity for solutions of the elliptic Stokes problem with the Neumann or Dirichlet boundary conditions (see [LT, Proposition 4.2] for a proof, see also [G, GS, S]).

**Lemma 3** (Basic $H^2$ and $H^3$ Stokes estimates). *Let $(u,p)$ be a solution of the problem*

$$\begin{cases} -\Delta u + \nabla p = f & \text{in } \Omega_f \\ \operatorname{div} u = g, & \text{in } \Omega_f \\ u = h_1 & \text{on } \Gamma_f \\ \partial_3 u + p\, e_3 = h_2 & \text{on } \Gamma_c. \end{cases}$$

*Then*

$$\|u\|_2 + \|p\|_1 \lesssim \|f\| + \|g\|_1 + \|h_1\|_{H^{3/2}(\Gamma_f)} + \|h_2\|_{H^{1/2}(\Gamma_c)}. \tag{2.9}$$

*If the boundary condition on $\Gamma_c$ is replaced by $u = h_3$, then*

$$\|u\|_2 + \|\nabla p\| \lesssim \|f\| + \|g\|_1 + \|h_1\|_{H^{3/2}(\Gamma_f)} + \|h_3\|_{H^{3/2}(\Gamma_c)} \tag{2.10}$$

*and*

$$\|u\|_3 + \|\nabla p\|_1 \lesssim \|f\|_1 + \|g\|_2 + \|h_1\|_{H^{5/2}(\Gamma_f)} + \|h_3\|_{H^{5/2}(\Gamma_c)}. \tag{2.11}$$

We also recall that a standard elliptic estimate (see [Gri, Theorem 1.6.1.5] or [LM, Chapter 2], for example) gives that if $\psi$ solves

$$\begin{cases} -\Delta \psi = f & \text{in } \Omega_e \\ \psi = h_1 & \text{on } \Gamma_e \\ \psi = h_2 & \text{on } \Gamma_c, \end{cases}$$

then

$$\|\psi\|_2 \lesssim \|f\| + \|h_1\|_{H^{3/2}(\Gamma_e)} + \|h_2\|_{H^{3/2}(\Gamma_c)} \tag{2.12}$$

and

$$\|\psi\|_3 \lesssim \|f\|_1 + \|h_1\|_{H^{5/2}(\Gamma_e)} + \|h_2\|_{H^{5/2}(\Gamma_c)}. \tag{2.13}$$

We now use (2.12) and (2.13) for the wave equation and Lemma 3 for the Stokes system, together with (2.7), to derive estimates for a solution $(v, w, q, \eta, a)$ of the system (1.1)–(1.10). We assume that $t \geq 0$ is such that

$$\begin{aligned} \Big( \|v\|_{L^1((0,t);H^3)} + \|v\|_{L^\infty((0,t);H^3)} + \|\nabla q\|_{L^\infty((0,t);H^1)} \\ + \|v_t\|_{L^1((0,t);H^2)} + \|q_t\|_{L^1((0,t);H^1)} \Big)(1 + \|v\|_{L^1((0,t);H^3)}) \leq \gamma, \end{aligned} \tag{2.14}$$

where $\gamma > 0$ is a sufficiently small constant. We fix $\gamma$ at the end of this section; we only need $\gamma$ to be sufficiently small so that the absorption arguments hold in the following estimates.

We first note that the Piola identity (2.8) gives $\partial_k(a_{ki}\int_{\Omega_f} q) = 0$, so that (1.6) can be rewritten as

$$-\Delta v_i + \partial_i q = -\partial_t v_i - \partial_j((\delta_{jl} - a_{jl}a_{kl})\partial_k v_i) + \partial_k\left((\delta_{ki} - a_{ki})\left(q - \int_{\Omega_f} q\right)\right), \quad i = 1,2,3 \quad (2.15)$$
$$\operatorname{div} v = (\delta_{ki} - a_{ki})\partial_k v_i$$

in $\Omega_f \times (0,T)$. Applying (2.11) to (2.15) gives

$$\|v\|_3 + \|\nabla q\|_1 \lesssim \|v_t\|_1 + \|\nabla((I-a)\nabla v)\|_1 + \|v\|_{H^{5/2}(\Gamma_c)} + \|\nabla((I-aa^T)\nabla v)\|_1$$
$$+ \left\|\nabla\left((I-a)\left(q - \int_{\Omega_f} q\right)\right)\right\|_1$$
$$\lesssim \|v_t\|_1 + \|I-a\|_2\|v\|_3 + \|v\|_{H^{5/2}(\Gamma_c)} + \|I-aa^T\|_2\left(\|v\|_3 + \left\|q - \int_{\Omega_f} q\right\|_2\right)$$
$$\lesssim \|v_t\|_1 + \|v\|_{H^{5/2}(\Gamma_c)} + \gamma\left(\|v\|_3 + \|\nabla q\|_1\right),$$

where we used (2.7) and the Poincaré inequality in the last line. Absorbing the last term by the left-hand side we obtain

$$\|v\|_3 + \|\nabla q\|_1 \lesssim \|v_t\|_1 + \|v\|_{H^{5/2}(\Gamma_c)} \lesssim \|v_t\|_1 + \|v\|_{H^{1/2}(\Gamma_c)} + \|\partial''v\|_{H^{1/2}(\Gamma_c)} \lesssim \|v_t\|_1 + \|w_t\|_1 + \|\partial'v\|_2. \quad (2.16)$$

In a similar way, but using (2.10) instead of (2.11), we obtain an $H^2$ estimate,

$$\|v\|_2 + \|\nabla q\| \lesssim \|v_t\| + \|v\|_{H^{3/2}(\Gamma_c)} \lesssim \|v_t\| + \|v\|_{H^{1/2}(\Gamma_c)} + \|\partial'v\|_{H^{1/2}(\Gamma_c)} \lesssim \|v_t\| + \|w_t\|_1 + \|\partial'w_t\|_1, \quad (2.17)$$

where we also used the continuity of the displacement across $\Gamma_c$, (1.8), and a trace estimate in the last inequalities of (2.16) and (2.17).

We now observe that $\partial' v$ satisfies the Stokes equation

$$-\Delta \partial' v_i + \partial_i \partial' q = -\partial' \partial_t v_i - \partial' \partial_j((\delta_{jl} - a_{jl}a_{kl})\partial_k v_i) + \partial' \partial_k\left((\delta_{ki} - a_{ki})\left(q - \int_{\Omega_f} q\right)\right)$$
$$\operatorname{div} \partial' v = \partial'((\delta_{ki} - a_{ki})\partial_k v_i)$$

(by taking $\partial'$ of (2.15)), with the boundary condition

$$\partial_3 \partial' v_i + \delta_{3i}\partial' q = \partial_3 \partial' w_i + \partial'((\delta_{3j} - a_{3l}a_{kl})\partial_k v_i) + \partial'((\delta_{3i} - a_{3i})q).$$

Thus, applying (2.10), we obtain

$$\|\partial' v\|_2 + \|\partial' q\|_1 \lesssim \|\partial' v_t\| + \|(I-a)\nabla v\|_2 + \|(I-aa^T)\nabla v\|_2$$
$$+ \left\|(I-a)\left(q - \int_{\Omega_f} q\right)\right\|_2 + \|\partial' v\|_{H^{3/2}(\Gamma_c)} \quad (2.18)$$
$$\lesssim \|\partial' v_t\| + \gamma\left(\|v\|_3 + \|\nabla q\|_1\right) + \|\partial''v\|_1 + \|\partial' v\|_1,$$

where we also used $\int_{\Omega_f} \partial' q = 0$ to write $\|\partial' q\|_1 \lesssim \|\nabla \partial' q\|$ on the left-hand side.

Alternatively, we may estimate $\|\partial' v\|_2 + \|\nabla \partial' q\|$ by (2.9), obtaining

$$\|\partial' v\|_2 + \|\partial' q\|_1 \lesssim \|\partial' v_t\| + \|(I - a)\nabla v\|_2 + \|(I - aa^T)\nabla v\|_2 + \left\|(I - a)\left(q - \int_{\Omega_f} q\right)\right\|_2 + \|\partial_3 \partial' w\|_{H^{1/2}(\Gamma_c)}$$
$$+ \|\partial'((I - aa^T)\nabla v)\|_{H^{1/2}(\Gamma_c)} + \|\partial' a_{3i} q\|_{H^{1/2}(\Gamma_c)} + \|(I - a)\partial' q\|_{H^{1/2}(\Gamma_c)}$$
$$\lesssim \|\partial' v_t\| + \gamma(\|v\|_3 + \|\nabla q\|_1) + \|\partial' w\|_2 + \|\partial' a_{3i}\|_{H^{1/2}(\Gamma_c)} \|q\|_{H^{3/2}(\Gamma_c)},$$
(2.19)

where, in the last inequality, we used the trace estimates to get

$$\|\partial'((I - aa^T)\nabla v)\|_{H^{1/2}(\Gamma_c)} + \|(I - a)\partial' q\|_{H^{1/2}(\Gamma_c)} \lesssim \|(I - aa^T)\nabla v\|_2 + \|(I - a)\nabla q\|_1 \lesssim \gamma(\|v\|_3 + \|\nabla q\|_1).$$

For the last term in (2.19), we observe that, by (2.3), we may write

$$\partial' a_{3i} = \frac{1}{2} \epsilon_{3mn} \epsilon_{ikl} \partial' \partial_m \eta_k \partial_n \eta_l + \frac{1}{2} \epsilon_{3mn} \epsilon_{ikl} \partial_m \eta_k \partial' \partial_n \eta_l,$$

and since $\epsilon_{3mn} \neq 0$ only if $m, n \in \{1, 2\}$, all the derivatives on the right-hand side are in the horizontal variables $x_1, x_2$. Thus recalling that $\nabla \eta = I + \int_0^t \nabla v$ (see (2.1)) and $w = \int_0^t w_t = \int_0^t v$ on $\Gamma_c$ we see that

$$\partial'' \eta = \int_0^t \partial'' v = \int_0^t \partial'' w_t = \partial'' w$$

on $\Gamma_c$, and consequently

$$\|\partial' a_{3i}\|_{H^{1/2}(\Gamma_c)} \lesssim \|\partial'' w\|_{H^{1/2}(\Gamma_c)} \|\nabla \eta\|_{H^{3/2}(\Gamma_c)} \lesssim \|\partial'' w\|_1 (1 + \|I - \nabla \eta\|_2) \lesssim \|\partial'' w\|_1, \qquad (2.20)$$

where we used (2.7) in the last inequality. Substituting (2.20) into (2.19), we get

$$\|\partial' v\|_2 + \|\partial' q\|_1 \lesssim \|\partial' v_t\| + \gamma(\|v\|_3 + \|\nabla q\|_1) + \|\partial' w\|_2 (1 + \|q\|_2). \qquad (2.21)$$

Applying (2.21) in (2.16) and absorbing the terms involving $\gamma$ we obtain

$$\|v\|_3 + \|\nabla q\|_1 \lesssim \|v_t\|_1 + \|w_t\|_1 + \|\partial' w\|_2 (1 + \|q\|_2). \qquad (2.22)$$

As for the pressure function, we note that $\|q\|$ can be estimated with a cumulative dependence on $v$ from the time 0 until a given time $t$. Namely,

$$\|q\| \lesssim \|\nabla q\| + \|q\|_{L^2(\Gamma_c)} \lesssim \|\nabla q\| + \|q(0)\|_{L^2(\Gamma_c)} + \int_0^t \|q_t\|_{L^2(\Gamma_c)}$$
$$\lesssim \|\nabla q\| + \|\partial_3 v(0)\|_{L^2(\Gamma_c)} + \int_0^t \|q_t\|_1 \lesssim \|\nabla q\| + \|v(0)\|_2 + \int_0^t (\|v_t\|_2 + \|q_t\|_1) \lesssim \gamma.$$
(2.23)

Thus, recalling (2.14), we have

$$\|q\|_2 \lesssim \gamma. \qquad (2.24)$$

For $w$, the elliptic estimate (2.13) applied to the Poisson equation $\Delta \overline{w} = \Delta w$, which holds due to $\partial' w = \partial' \overline{w}$ and $\partial_{33} w = \partial_{33} \overline{w}$, gives

$$\|\overline{w}\|_3 \lesssim \|\Delta w\|_1 + \|\overline{w}\|_{H^{5/2}(\Gamma_c)} \lesssim \|w_{tt}\|_1 + \|w_t\|_1 + \|\partial'' \overline{w}\|_1 = \|w_{tt}\|_1 + \|w_t\|_1 + \|\partial'' w\|_1, \qquad (2.25)$$

where in the second inequality we used (1.1) and the fact that both $\int_{\Gamma_c} \overline{w}$ and $\int_{\Gamma_c} \partial' \overline{w}$ vanish (recall (1.18)). Moreover, the elliptic estimate (2.12) applied to $\Delta w_t = w_{ttt} + \alpha w_{tt}$ gives

$$\|w_t\|_2 \lesssim \|w_{ttt}\| + \|w_{tt}\| + \|w_t\|_{H^{3/2}(\Gamma_c)} \lesssim \|w_{ttt}\| + \|w_{tt}\| + \|v\|_2. \tag{2.26}$$

As for $v_t$, we see that it satisfies

$$-\Delta \partial_t v_i + \partial_i q_t = -\partial_{tt} v_i - \partial_t \partial_j ((\delta_{jk} - a_{jl} a_{kl}) \partial_k v_i) + \partial_t ((\delta_{ji} - a_{ji}) \partial_j q)$$
$$\operatorname{div} v_t = \partial_t \partial_k ((\delta_{ki} - a_{ki}) v_i)$$

with the boundary conditions

$$\partial_3 \partial_t v_i + \delta_{3i} \partial_t q = \partial_3 \partial_t w_i + \partial_t ((\delta_{3j} - a_{3l} a_{kl}) \partial_k v_i) + \partial_t ((\delta_{3i} - a_{3i}) q).$$

Thus (2.9) gives

$$\begin{aligned}\|v_t\|_2 + \|q_t\|_1 &\lesssim \|v_{tt}\| + \|\partial_t \nabla((I-a)\nabla v)\|_1 + \|\partial_t \nabla((I-aa^T)\nabla v)\| + \|\partial_t((I-a)\nabla q)\| \\ &\quad + \|\partial_3 w_t\|_{H^{1/2}(\Gamma_c)} + \|\partial_t((I-aa^T)\nabla v)\|_{H^{1/2}(\Gamma_c)} + \|\partial_t((I-a)q)\|_{H^{1/2}(\Gamma_c)} \\ &\lesssim \|v_{tt}\| + \|a_t\|_2 \|v\|_2 + \|I-a\|_2 \|v_t\|_2 + \|I-aa^T\|_2 \|v_t\|_2 + \|\partial_t(aa^T)\|_2 \|v\|_2 \\ &\quad + \|a_t\|_2 \|\nabla q\| + \|I-a\|_2 \|q_t\|_1 + \|w_t\|_2 \\ &\lesssim \|v_{tt}\| + \gamma(\|v_t\|_2 + \|q_t\|_1) + \|v\|_3(\|v\|_2 + \|\nabla q\|) + \|w_{ttt}\| + \|w_{tt}\| + \|v\|_2,\end{aligned} \tag{2.27}$$

where in the last inequality we used (2.26) as well as (2.2), (2.3), (2.4), and (2.7) to write

$$\|\partial_t(aa^T)\|_2 \lesssim \|a_t\|_2 \|a\|_2 \lesssim \|\nabla \eta\|_2 \|v\|_3 (1 + \|I-a\|_2) \lesssim \|v\|_3$$

and similarly for $\|a_t\|_2$. Thus, absorbing the terms with $\gamma$ by the left-hand side, we obtain

$$\|v_t\|_2 + \|q_t\|_1 \lesssim \|v_{tt}\| + \gamma \|v\|_3 + \|w_{ttt}\| + \|w_{tt}\| + \|v\|_2, \tag{2.28}$$

where we also used $\|v\|_2 + \|\nabla q\| \lesssim \gamma$ on the third term in the last line of (2.27).

Thus using the pressure estimate (2.24), the $H^2$ estimate on $\partial' v$ (2.21), the $H^2$ estimate (2.17) we can combine the $H^3$ estimate (2.22), the $H^2$ estimate of $v_t$ (2.28), and the $H^3$ bound for $\overline{w}$, (2.25), into

$$\|v\|_3 + \|\nabla q\|_1 + \|v_t\|_2 + \|q_t\|_1 + \|\overline{w}\|_3 \lesssim \|v_{tt}\| + \|\partial' v_t\| + \|v_t\| + \|w_{ttt}\| + \|w_{tt}\|_1 + \|\partial' w_t\|_1 + \|w_t\|_1 + \|\partial'' w\|_1. \tag{2.29}$$

We now fix $\gamma > 0$ to be sufficiently small so that the above absorption arguments are valid. This last estimate (2.29) and (2.26) is very useful for us, as the right-hand sides involve only the terms appearing in (1.13) for $S \ne \mathrm{id}$. Thus we are able (in (3.15) below) to control all the relevant Sobolev norms of our model (1.1)–(1.9) by the total energy of the system $X(t)$.

# 3 Proof of the main result

Here we prove our main result, Theorem 1. We shall consider

$$S \in \{\partial', \partial_t, \partial'\partial_t, \partial'', \partial_{tt}\},$$

and for each $S$ derive an energy estimate, stated in Section 3.1. In Section 3.2 we then use such individual energy estimates to define the total energy of the system and estimate it in the form of an ODE-type inequality. Then, using the ODE-type Lemma 4, which we prove in Section 3.3, we obtain the required a priori bound for the total energy.

## 3.1 Energy estimates

Here we derive an energy inequality for each $S \in \{\partial', \partial_t, \partial'\partial_t, \partial'', \partial_{tt}\}$. First we apply $S$ to (1.6) to get

$$S\partial_t v_i - \partial_j(a_{jl}a_{kl}\partial_k Sv_i) + a_{ki}S\partial_k q$$
$$= \left(S\partial_j(a_{jl}a_{kl}\partial_k v_i) - \partial_j(a_{jl}a_{kl}\partial_k Sv_i)\right) - \left(S(a_{ki}\partial_k q) - a_{ki}\partial_k Sq\right) \text{ in } \Omega_\mathrm{f} \times (0,T), \qquad i = 1,2,3,$$

and we test it with $Sv_i$. We also apply $S$ to (1.1) and test it with $\partial_t Sw_i$. Summing the resulting equations we get

$$\frac{1}{2}\frac{\mathrm{d}}{\mathrm{d}t}\Big(\|Sv\|_{L^2}^2 + \|Sw_t\|_{L^2}^2 + \|S\nabla w\|_{L^2}^2\Big) + \|\nabla Sv\|_{L^2}^2 + \alpha\|Sw_t\|_{L^2}^2$$
$$= \int_{\Omega_\mathrm{f}} (S\partial_j(a_{jl}a_{kl}\partial_k v_i) - \partial_j(a_{jl}a_{kl}\partial_k Sv_i))Sv_i - \int_{\Omega_\mathrm{f}} (S(a_{ki}\partial_k q) - a_{ki}\partial_k Sq)Sv_i$$
$$+ \int_{\Omega_\mathrm{f}} (\delta_{jk} - a_{jl}a_{kl})\partial_k Sv_i\partial_j Sv_i + \underbrace{\int_{\Omega_\mathrm{f}} Sq\partial_k(a_{ki}Sv_i)}_{\overset{(1.7)\&(2.8)}{=}\int_{\Omega_\mathrm{f}} Sq(a_{ki}S\partial_k v_i - S(a_{ki}\partial_k v_i))} \quad . \qquad (3.1)$$

We now derive the lower level energy estimate. Namely, we apply $S$ to (1.1) and test it with $Sw$ to obtain

$$\frac{\mathrm{d}}{\mathrm{d}t}\int_{\Omega_\mathrm{e}} Sw_t \cdot Sw - \int_{\Omega_\mathrm{e}} |Sw_t|^2 + \int_{\Omega_\mathrm{e}} \partial_j Sw_i \partial_j Sw_i + \frac{\alpha}{2}\frac{\mathrm{d}}{\mathrm{d}t}\int_{\Omega_\mathrm{e}} |Sw|^2 - \int_{\Gamma_\mathrm{c}} \partial_3 Sw_i Sw_i \mathrm{d}\sigma(x) = 0,$$

or, equivalently,

$$\frac{\mathrm{d}}{\mathrm{d}t}\left(\frac{\alpha}{2}\|Sw\|^2 + \int_{\Omega_\mathrm{e}} \partial_t Sw \cdot Sw\right) + \|\nabla Sw\|^2 = \|Sw_t\|^2 + \int_{\Gamma_\mathrm{c}} \partial_3 Sw_i Sw_i\, d\sigma(x). \qquad (3.2)$$

Given $\tau \geq 0$ and $t > \tau$ we now consider the test function

$$\phi(t) := S\eta(t) - S\eta(\tau) + \widetilde{Sw}(\tau) \quad \text{on } \Omega_\mathrm{f}, \qquad (3.3)$$

where $\widetilde{f}$ denotes an extension of $f$, defined on $\Omega_\mathrm{e}$, to the fluid domain $\Omega_\mathrm{f}$, such that

$$\widetilde{f}\,\big|_{\Gamma_\mathrm{c}} = f\,\big|_{\Gamma_\mathrm{c}}, \qquad \widetilde{f}\,\big|_{\Gamma_\mathrm{f}} = 0, \qquad \text{and} \qquad \|\widetilde{f}\|_1 \lesssim \|f\|_1. \qquad (3.4)$$

Note that integrating (1.8) on the time interval $(\tau, t)$ gives

$$\eta(t) - \eta(\tau) + w(\tau) = w(t) \quad \text{on } \Gamma_c,$$

and so applying $S$ gives that

$$\phi(t)\big|_{\Gamma_c} = Sw(t)\big|_{\Gamma_c} = \widetilde{Sw}(t)|_{\Gamma_c}. \tag{3.5}$$

Applying $S$ to the velocity equation (1.6) and testing it with $\phi$, we obtain

$$\int_{\Omega_f} S\partial_t v_i \phi_i - \int_{\Omega_f} \partial_j S(a_{jl} a_{kl} \partial_k v_i) \phi_i + \int_{\Omega_f} \partial_k S(a_{ki} q) \phi_i = 0,$$

from where, by $\partial_t \phi_i = S\partial_t \eta_i = Sv_i$,

$$\begin{aligned}
&\frac{\mathrm{d}}{\mathrm{d}t} \int_{\Omega_f} Sv_i \phi_i - \int_{\Omega_f} |Sv|^2 + \int_{\Omega_f} S(a_{jl} a_{kl} \partial_k v_i) \partial_j \phi_i \\
&+ \int_{\Gamma_c} S(a_{jl} a_{kl} \partial_k v_i) \phi_i N_j \, d\sigma(x) - \int_{\Omega_f} S(a_{ki} q) \partial_k \phi_i - \int_{\Gamma_c} S(a_{ki} q) \phi_i N_k \, d\sigma(x) = 0.
\end{aligned} \tag{3.6}$$

We rewrite the third term as

$$\begin{aligned}
\int_{\Omega_f} S(a_{jl} a_{kl} \partial_k v_i) \partial_j \phi_i &= \int_{\Omega_f} \nabla Sv : \nabla \phi + \int_{\Omega_f} S((a_{jl} a_{kl} - \delta_{jk}) \partial_k v_i) \partial_j \phi_i \\
&= \int_{\Omega_f} \nabla S \partial_t (\eta - \eta(\tau)) : \nabla S(\eta - \eta(\tau)) + \int_{\Omega_f} \nabla S \partial_t (\eta - \eta(\tau)) : \nabla \widetilde{Sw}(\tau) \\
&\quad + \int_{\Omega_f} S((a_{jl} a_{kl} - \delta_{jk}) \partial_k v_i) \partial_j \phi_i \\
&= \frac{1}{2} \frac{\mathrm{d}}{\mathrm{d}t} \|\nabla S(\eta - \eta(\tau))\|^2 + \int_{\Omega_f} \nabla Sv : \nabla \widetilde{Sw}(\tau) + \int_{\Omega_f} S((a_{jl} a_{kl} - \delta_{jk}) \partial_k v_i) \partial_j \phi_i,
\end{aligned} \tag{3.7}$$

where $A : B = A_{ij} B_{ij}$, and the fifth term in (3.6) as

$$-\int_{\Omega_f} S(a_{ki} q) \partial_k \phi_i \, dx = -\int_{\Omega_f} Sq \,\mathrm{div}\, \phi + \int_{\Omega_f} S((\delta_{ki} - a_{ki}) q) \partial_k \phi_i. \tag{3.8}$$

Applying (3.7) and (3.8) in (3.6) and adding the resulting equation to (3.2), we observe that the integrals over $\Gamma_c$ cancel due to (3.5) and the boundary condition (1.9). We obtain

$$\begin{aligned}
&\frac{\mathrm{d}}{\mathrm{d}t} \left( \frac{\alpha}{2} \|Sw\|^2 + \int_{\Omega_e} Sw_t \cdot Sw + \frac{1}{2} \|\nabla S(\eta - \eta(\tau))\|^2 + \int_{\Omega_f} Sv \cdot \phi \right) + \|\nabla Sw\|^2 \\
&= \|Sw_t\|^2 + \|Sv\|^2 - \int_{\Omega_f} \nabla Sv : \nabla \widetilde{Sw}(\tau) + \int_{\Omega_f} S((\delta_{jk} - a_{jl} a_{kl}) \partial_k v_i) \partial_j \phi_i \\
&\quad + \int_{\Omega_f} Sq \,\mathrm{div}\, \phi - \int_{\Omega_f} S((\delta_{ki} - a_{ki}) q) \partial_k \phi_i.
\end{aligned} \tag{3.9}$$

We now multiply (3.9) by $\lambda \in (0,1]$ (to be fixed later) and add to (3.1) to obtain

$$\frac{\mathrm{d}}{\mathrm{d}t} E_S(t, \tau) + D_S(t) \leq L_S(t, \tau) + N_S(t, \tau) + C_S(t), \tag{3.10}$$

for each $S \in \{\partial', \partial_t, \partial'\partial_t, \partial'', \partial_{tt}\}$, where the pointwise energy $E_S$, the dissipation energy $D_S$, the linear part $L_S$, the nonlinear part $N_S$, and the commutator part $C_S$ are defined by

$$\begin{aligned}
E_S(t,\tau) &:= \frac{1}{2}\|Sv(t)\|^2 + \frac{1}{2}\|Sw_t(t)\|^2 + \frac{1}{2}\|\nabla Sw(t)\|^2 \\
&\quad + \frac{\lambda\alpha}{2}\|Sw\|^2 + \frac{\lambda}{2}\|\nabla S(\eta-\eta(\tau))\|^2 + \lambda\int_{\Omega_{\mathrm{f}}} Sv\cdot\phi + \lambda\int_{\Omega_{\mathrm{e}}} Sw_t \cdot Sw, \\
D_S(t) &:= \frac{1}{2}\|\nabla Sv\|^2 + \frac{1}{C}\|Sv\|^2 + (\alpha-\lambda)\|Sw_t\|^2 + \frac{\lambda}{2}\|\nabla Sw\|^2 + \frac{\lambda}{C}\|Sw\|^2, \\
L_S(t,\tau) &:= -\lambda\int_{\Omega_{\mathrm{f}}} \nabla Sv : \nabla\widetilde{Sw}(\tau) + \lambda\int_{\Omega_{\mathrm{f}}} Sq\,\mathrm{div}\,\phi =: L_{S,1}(t,\tau) + L_{S,2}(t,\tau), \\
N_S(t,\tau) &:= \lambda\int_{\Omega_{\mathrm{f}}} S((\delta_{jk} - a_{jl}a_{kl})\partial_k v_i)\partial_j\phi_i - \lambda\int_{\Omega_{\mathrm{f}}} S((\delta_{ki} - a_{ki})q)\partial_k\phi_i + \int_{\Omega_{\mathrm{f}}} (\delta_{jk} - a_{jl}a_{kl})\partial_k Sv_i \partial_j Sv_i, \\
C_S(t) &:= \int_{\Omega_{\mathrm{f}}} (S\partial_j(a_{jl}a_{kl}\partial_k v_i) - \partial_j(a_{jl}a_{kl}\partial_k Sv_i))Sv_i \\
&\quad - \int_{\Omega_{\mathrm{f}}} (S(a_{ki}\partial_k q) - a_{ki}\partial_k Sq)Sv_i + \int_{\Omega_{\mathrm{f}}} Sq\,(a_{ki}S\partial_k v_i - S(a_{ki}\partial_k v_i)),
\end{aligned} \tag{3.11}$$

and where $\phi$, which depends on $\tau$, $t$ and $S$, is defined in (3.3) and $C \geq 1$ is a constant; note that we used the Poincaré inequalities $\|Sv\|^2 \lesssim \|\nabla Sv\|^2$ and $\|Sw\|^2 \lesssim \|\nabla Sw\|^2$ to include the terms $\|Sv\|^2$ and $\|Sw\|^2$ in $D_S$.

By applying Young's inequality $\lambda ab \leq a^2/4 + \lambda^2 b^2$ to the last two terms of $E_S$ and using the Poincaré inequality $\|S(\eta-\eta(\tau))\| \lesssim \|\nabla S(\eta-\eta(\tau))\|$ we see that $E_S$ may be bounded from above and below as

$$E_S(t,\tau) \sim \|Sv(t)\|^2 + \|Sw_t(t)\|^2 + \|Sw(t)\|_1^2 + \lambda\|\nabla S(\eta-\eta(\tau))\|^2 + \lambda\int_{\Omega_{\mathrm{f}}} Sv\cdot\widetilde{Sw}(\tau), \tag{3.12}$$

as long as $\lambda > 0$ is chosen sufficiently small, where $a \sim_\alpha b$ means $a \lesssim_\alpha b$ and $b \lesssim_\alpha a$.

## 3.2 The total energy of the system

In this section we combine the energy estimates (3.10) from the previous section to define the total energy $X(t)$ of the system (1.1)–(1.10) and derive an a priori estimate for it.

Before defining $X$, we note that we aim to obtain an ODE-type estimate of the form

$$X(t) + \lambda\int_\tau^t X \lesssim \left(1 + \lambda^2(t-\tau)\right)X(\tau) + \lambda^2\int_\tau^t X + \text{(small terms)}, \tag{3.13}$$

for all $\lambda > 0$, $\tau \geq 0$ and $t \geq \tau$, where the small terms are at least cubic in $X^{1/2}$. Note that (3.13) leaves the choice of $\lambda > 0$ free. This way, we are able to choose $\lambda$ small enough so that the factor of "$t-\tau$" and the linear term $\int_\tau^t X$ on the right-hand side are negligible in proving an exponential decay of $X(t)$ for a sufficiently small $X(0)$. We make this precise in Lemma 4 below.

The question therefore becomes which terms appearing in the energy estimates (3.10), where $S \in \{\partial', \partial_t, \partial'\partial_t, \partial'', \partial_{tt}\}$, should be included in the total energy $X(t)$. To this end we observe that (3.13) requires that, given $\tau \geq 0$, we must be able to control each such term for every $t \geq \tau$ as well as its integral from $\tau$ to $t$. This can be ensured by choosing only terms that appear in both $E_S$ and $D_S$, for some $S \in \{\partial', \partial_t, \partial'\partial_t, \partial'', \partial_{tt}\}$. Thus we let

$$X(t) \coloneqq \sum_{S \in \{\partial', \partial_t, \partial'\partial_t, \partial'', \partial_{tt}\}} \left( \|Sv(t)\|^2 + \|Sw_t(t)\|^2 + \|Sw(t)\|_1^2 \right). \tag{3.14}$$

With this definition of $X$, we can control all the quantities appearing in the energy estimates (3.10) by $X$, that is, there exists a universal constant $\overline{C} \geq 1$ such that

$$\|v\|_3^2 + \|\nabla q\|_1^2 + \|v_t\|_2^2 + \|q_t\|_1^2 + \|\overline{w}\|_3^2 + \|w_t\|_2^2 + \|w_{tt}\|_1^2 + \|\partial' w\|_2^2 \leq \overline{C} X, \tag{3.15}$$

for all times $t \geq 0$ such that

$$h(t) \coloneqq \sup_{[0,t)} X^{1/2} + \int_0^t X^{1/2} \leq \frac{\gamma}{10\overline{C}}. \tag{3.16}$$

Recall (2.14) and that $\gamma$ is the constant fixed in Section 2. The estimate (3.15) is essentially a consequence of (2.26) and (2.29), but it requires a comment. In fact, in order to see (3.15)–(3.16) we first note that (3.15) is satisfied with some $\overline{C}$ at $t = 0$, as then $\eta = x$, $a = I$, $w = 0$, and so the Stokes and elliptic estimates (2.9)–(2.13) apply directly, without using (2.14). Taking $\overline{C}$ larger we deduce, by a continuity argument, that (2.14) holds also for some $t > 0$. Secondly, supposing that $\overline{C}$ is larger than the constant resulting from combining (2.26) and (2.29), we see that, as long as (3.15) remains valid, (3.16) implies (2.14). Thus (2.26) and (2.29) are valid and imply that (3.15) continues to hold, which lets us use a continuity argument to deduce (3.15)–(3.16).

We note that (3.15) involves all variables of the system (1.1)–(1.10), except for the lowest order terms $\|q\|$ and $\|w\|_1$. As mentioned in the introduction (below Theorem 2), we cannot expect these quantities to decay to 0. Instead they can, roughly speaking, accumulate from time 0. In fact, the purpose of introducing the function $h(t)$ above is to keep track of such quantities. For example, in the same way as in (2.24) we obtain

$$\|q(t)\| + \|w(t)\|_1 \lesssim h(t) \tag{3.17}$$

given (3.16) holds.

Let us comment further on our choice of $X$ in (3.14) and on the energy estimates (3.10) that motivate our strategy for proving Theorem 1. We first note that some terms appear in $E_S$, but not in $D_S$, and vice versa. For example, $\|\nabla S(\eta - \eta(\tau))\|$ appears in $E_S$ for each $S$ (recall (3.12)), but we shall not use this term. In fact, if $S \in \{\partial', \partial''\}$ then term can simply be ignored, as it vanishes at time $\tau$. If, on the other hand, $S$ involves time derivative, then we observe that $\|\nabla Sv\|$ is included in the dissipation $D_S$ without a factor $\lambda$, which can be used to neglect this term (see (3.23) below). We also note that the term $\|Sw\|$, which is included in $E_S$ with the (small) weight $\lambda$, but this can be neglected as well due to the Poincaré inequality, since $\|\nabla Sw\|$ is also included in $E_S$ with weight 1.

We also observe that $\|\nabla Sw\|$ is included in $D_S$ with the weight $\lambda$, and this is the reason for the appearance of a factor of $\lambda$ on the left-hand side of (3.13). Consequently this influences the exponential decay rate of $X(t)$ as $t \to \infty$ (see Lemma 4 for details). Moreover the appearance of this factor also requires that all linear terms on the right-hand side of (3.13) with a factor of $t - \tau$ are also accompanied by $\lambda^2$ (or, more generally, $\lambda^b$ for some $b > 1$), in order to make sure that they are negligible compared to the left-hand side of (3.13). Indeed, choosing $\lambda > 0$ sufficiently small then enables one to obtain exponential decay (see Lemma 4 for details).

We further observe that $\|\nabla Sv\|^2$ appears in $D_S$, but not in $E_S$, for each $S$. This can be thought of as the main dissipation term, and we make full use of it in estimating some linear parts $L_S$, as well as $N_{\partial_{tt}}$. In fact, we show in Step 1 below that, given (3.16) holds,

$$L_{\partial'} + L_{\partial''} \lesssim \delta \left( \|v_t\|_1^2 + \|\partial'' v\|_1^2 + \|\partial' v\|_1^2 \right) + C_\delta \lambda^2 X(\tau) + hO(X), \tag{3.18}$$

and in Step 2 below that

$$\int_\tau^t (N_{\partial_{tt}} + L_{\partial_{tt}}) \lesssim \delta \int_\tau^t \|\nabla \partial_{tt} v\|^2 + C_\delta hO(X), \tag{3.19}$$

for all $\delta > 0$, where

$$\begin{aligned}&O(X^{k/2}) \text{ denotes any finite sum of products of any power of } t - \tau, \text{ any power of } \lambda \text{ and} \\ &\text{at least } k \text{ factors of the form } X(t)^{1/2}, X(\tau)^{1/2} \text{ or } \int_\tau^t X^{1/2},\end{aligned} \tag{3.20}$$

for $k \geq 0$. Furthermore, we show in Step 3 below that, given (3.16) holds, all the other terms appearing on the right hand sides of (3.10) for all $S$ are negligible, namely that

$$\int_\tau^t \left( L_{\partial_t} + L_{\partial' \partial_t} + \sum_{S \in \{\partial', \partial_t, \partial' \partial_t, \partial''\}} N_S + \sum_{S \in \{\partial', \partial_t, \partial' \partial_t, \partial'', \partial_{tt}\}} C_S \right) \lesssim hO(X). \tag{3.21}$$

With the estimates (3.18)–(3.21) in hand, we can conclude the proof of Theorem 1 using the strategy described above.

*Proof of Theorem 1.* We first sum the estimates (3.10) for all $S \in \{\partial', \partial_t, \partial' \partial_t, \partial'', \partial_{tt}\}$ and recall (3.12)

to obtain

$$\begin{aligned}
X(t) + &\sum_{S \in \{\partial_t, \partial', \partial'\partial_t, \partial_{tt}, \partial''\}} \int_\tau^t \left( \|Sv\|_1^2 + \|Sw_t\|^2 + \lambda \|Sw\|_1^2 \right) \\
&\qquad + \lambda \|\nabla \partial'(\eta(t) - \eta(\tau))\|^2 + \lambda \|\nabla \partial''(\eta(t) - \eta(\tau))\|^2 \\
\lesssim\ & X(\tau) + \lambda \left( \left( \|\nabla v(t)\|^2 + \|\nabla \partial' v(t)\|^2 + \|\nabla v_t(t)\|^2 \right) - \left( \|\nabla v(\tau)\|^2 + \|\nabla \partial' v(\tau)\|^2 + \|\nabla v_t(\tau)\|^2 \right) \right) \\
&\qquad + \sum_{S \in \{\partial_t, \partial', \partial'\partial_t, \partial_{tt}, \partial''\}} \int_\tau^t (L_S + N_S + C_S) \\
\lesssim\ & X(\tau) + \delta \int_\tau^t \left( \|v_t\|_1^2 + \|\partial' v_t\|_1^2 + \|v_{tt}\|_1^2 + \|\partial'' v\|_1^2 + \|\partial' v\|_1^2 \right) \\
&\qquad + C_\delta \lambda^2 \left( (t-\tau) X(\tau) + \int_\tau^t X \right) + h O(X),
\end{aligned} \tag{3.22}$$

where, in the last inequality, we have used (3.18)–(3.21), as well as estimated the second term following the first inequality by

$$2\lambda \int_\tau^t \left( \|\nabla v\| \|\nabla v_t\| + \|\nabla \partial' v\| \|\nabla \partial' v_t\| + \|\nabla v_t\| \|\nabla v_{tt}\| \right)$$
$$\lesssim \delta \int_\tau^t \left( \|\nabla v_t\|^2 + \|\nabla \partial' v_t\|^2 + \|\nabla v_{tt}\|^2 \right) + C_\delta \lambda^2 \int_\tau^t X,$$

due to the Fundamental Theorem of Calculus and (3.15). Note that we have also absorbed the terms $\lambda \int_{\Omega_\mathrm{f}} Sv \cdot \widetilde{Sw}(\tau)$, included in $E_S$ (recall (3.12)), by estimating

$$\lambda \int_{\Omega_\mathrm{f}} (Sv(t) - Sv(\tau)) \cdot \widetilde{Sw}(\tau) \leq \lambda^2 X(t) + C X(\tau)$$

and absorbing the first of the resulting terms by the left-hand side.

We now neglect the last two terms on the far left side of (3.22), and fix a sufficiently small $\delta > 0$ so that the second term on the far right side can be absorbed by the far left side. Then noting that the second term on the left-hand side is larger than or equal to $\lambda \int_\tau^t X$ we obtain

$$X(t) + \lambda \int_\tau^t X \leq C \left( 1 + \lambda^2 (t-\tau) \right) X(\tau) + C\lambda^2 \int_\tau^t X + h O(X), \tag{3.23}$$

for all $\lambda \in (0, 1]$, $\tau \geq 0$, and $t \geq \tau$.

As mentioned above, we emphasize that it is necessary that each term on the right-hand side of (3.23) that is linear in $X$ and includes a power of $(t-\tau)$ must also be accompanied by a factor of $\lambda^2$, as otherwise an exponential decay cannot be obtained (see Lemma 4 below). This in particular demonstrates the necessity of the factor of $\lambda^2$ (rather than merely $\lambda$) appearing at all terms that are linear in $X$ in (3.18). We now fix

$$\lambda := \frac{1}{500 C^2}, \tag{3.24}$$

where $C \geq 1$ is the constant in (3.23).

We also recall (3.16), which reminds us that (3.23) is only valid for $t$ such that $h(t) = \sup_{(0,t)} X^{1/2} + \int_0^t X^{1/2} \leq \gamma/10\overline{C}$. However, Lemma 4 shows that there exists $\varepsilon > 0$ such that if $X(0) \leq \varepsilon$ then both (3.16) and the estimate $X(t) \leq 30C \mathrm{e}^{-t/1000C^3}$ hold for all $t > 0$, as required. $\square$

It remains to prove (3.18), (3.19), and (3.21).

Step 1. We prove (3.18).

We first note that in this step we need to control $\|v\|_3 + \|\nabla q\|_1$ by a sum of all dissipation terms $D_S$. One may be inclined to simply use (3.15) since all the terms included in the definition (3.14) are also included in the dissipation terms. Such approach, however, would result in an additional factor of $\lambda^{-1}$, since terms $\|Sw\|^2$ and $\|\nabla Sw\|^2$ appear with the coefficient $\lambda$ in the definition (3.11) of $D_S$. Thus such control would be insufficient, as it would consequently give an ODE-type estimate on $X$ of the form of (3.13), but without one of the $\lambda$'s on the right-hand side, and so would make it impossible to apply Lemma 4 to conclude the proof.

Instead, we first write
$$\begin{aligned} L_{\partial'',1} + L_{\partial',1} &\lesssim \lambda \|\nabla \partial'' v\| \, \|\partial'' w(\tau)\|_1 + \lambda \|\nabla \partial' v\| \, \|\partial' w(\tau)\|_1 \\ &\lesssim \delta \left( \|\nabla \partial'' v\|^2 + \|\nabla \partial' v\|^2 \right) + C_\delta \lambda^2 \|\partial' w(\tau)\|_2^2 \lesssim \delta \left( \|\nabla \partial'' v\|^2 + \|\nabla \partial' v\|^2 \right) + C_\delta \lambda^2 X(\tau). \end{aligned}$$

For the second linear term, $L_{S,2}$, we first note that we may replace $\gamma$ by $h$ in (2.18) to obtain
$$\|\partial' v\|_2 + \|\partial' q\|_1 \lesssim \|\nabla v_t\| + \|\partial'' v\|_1 + \|\partial' v\|_1 + h \underbrace{(\|v\|_3 + \|\nabla q\|_1)}_{\lesssim X^{1/2}}. \tag{3.25}$$

We also note that applying the bound (3.15) in (2.2), (2.4), (2.5), and (2.6) we can replace $\gamma$ by $h$ in (2.7) to obtain
$$\begin{aligned} \|\nabla \eta\|_2, \|a\|_2 &\lesssim 1, \\ \|a_t\|_2, \|a_{tt}\|_1 &\lesssim X^{1/2}, \\ \|I - a\|_2, \|I - aa^T\|_2 &\lesssim h \end{aligned} \tag{3.26}$$

at each time instant $t$. Moreover, we observe that
$$\operatorname{div} \partial''(\eta - \eta(\tau)) = \int_\tau^t \partial'' (\delta_{ki} \partial_k v_i) = \int_\tau^t \partial'' ((\delta_{ki} - a_{ki}) \partial_k v_i),$$

and similarly when $\partial''$ is replaced by $\partial'$. Thus we obtain
$$\begin{aligned} L_{\partial'',2} + L_{\partial',2} &\lesssim \lambda \|\partial' q\|_1 (\|\operatorname{div} \partial''(\eta - \eta(\tau))\| + \|\operatorname{div} \partial'(\eta - \eta(\tau))\| + \|\partial'' \nabla w(\tau)\| + \|\nabla \partial' w(\tau)\|) \\ &\lesssim \lambda \|\partial' q\|_1 \left( \int_\tau^t \|(I-a)\nabla v\|_2 + X^{1/2}(\tau) \right) \\ &\lesssim \delta \left( \|v_t\|_1^2 + \|\partial'' v\|_1^2 + \|\partial' v\|_1^2 \right) + C_\delta \lambda^2 X(\tau) + hO(X), \end{aligned}$$

as required, where we used (3.4) in the first inequality and (3.15) in the second; we also applied (3.25) to bound $\|\partial' q\|_1$ and obtain the terms with $\delta$, as well as (3.26) and (3.15) to obtain the last term.

**Step 2.** We prove (3.19).

For $N_{\partial_{tt}}$ we first note that

$$\|\partial_t((I-a)q)\|_1 \lesssim \|a_t\|_2\|q\|_2 + \|I-a\|_2\|q_t\|_1 \lesssim hX^{1/2},$$
$$\|\partial_t((I-aa^T)\nabla v)\|_1 \lesssim \|a_t\|_2\|a\|_2\|v\|_3 + \|I-aa^T\|_2\|v_t\|_2 \lesssim hX^{1/2},$$
(3.27)

due to (3.15) and (3.26). Since

$$\int_\tau^t N_{\partial_{tt}} = \lambda \int_\tau^t \int_{\Omega_f} \partial_{tt}((\delta_{jk} - a_{jl}a_{kl})\partial_k v_i)\partial_j \partial_t v_i - \lambda \int_\tau^t \int_{\Omega_f} \partial_{tt}((\delta_{ki} - a_{ki})q)\partial_k \partial_t v_i$$
$$+ \int_\tau^t \int_{\Omega_f} (\delta_{jk} - a_{jl}a_{kl})\partial_k \partial_{tt} v_i \partial_j \partial_{tt} v_i,$$

we can integrate by parts in time in the first two integrals and then bound the resulting terms to obtain

$$\int_\tau^t N_{\partial_{tt}} \leq \lambda \int_\tau^t \|\nabla v_{tt}\| \left(\|\partial_t((I-aa^T)\nabla v)\| + \|\partial_t((I-a)q)\|\right) + \int_\tau^t \|I-aa^T\|_{L^\infty}\|\nabla \partial_{tt} v\|^2$$
$$+ \lambda \left[\int_{\Omega_f} \partial_t((\delta_{jk} - a_{jl}a_{kl})\partial_j v_i)\partial_k \partial_t v_i - \int_{\Omega_f} \partial_t((\delta_{ki} - a_{ki})q)\partial_k \partial_t v_i\right]_\tau^t$$
$$\lesssim \delta \int_\tau^t \|\nabla \partial_{tt} v\|^2 + C_\delta hO(X),$$

where we recalled the convention (3.20) applied Young's inequality and (3.27) in the last step.

For $L_{\partial_{tt}}$, we first note that

$$\|\partial_t((I-a)\nabla v)\|_1 \lesssim \|a_t\|_2\|v\|_3 + \|I-a\|_2\|v_t\|_2 \lesssim X + hX^{1/2},$$
$$\|\partial_{tt}((I-a)\nabla v)\| \lesssim \|I-a\|_2\|\nabla v_{tt}\| + \|a_t\|_2\|v_t\|_1 + \|a_{tt}\|\|v\|_3 \lesssim X + hX^{1/2},$$
(3.28)

as in (3.27), where we used (2.4) to estimate $\|a_{tt}\|_1 \lesssim O(X^{1/2})$. Therefore,

$$\int_\tau^t L_{\partial_{tt}} = \lambda \int_\tau^t \int_{\Omega_f} q_{tt} \operatorname{div} \partial_{tt} \eta = \lambda \int_\tau^t \int_{\Omega_f} q_{tt} \partial_t(\delta_{jk}\partial_j v_k) = \lambda \int_\tau^t \int_{\Omega_f} q_{tt} \partial_t((\delta_{jk} - a_{jk})\partial_j v_k)$$
$$= -\lambda \int_\tau^t \int_{\Omega_f} q_t \partial_{tt}((\delta_{jk} - a_{jk})\partial_j v_k) + \lambda \left[\int_{\Omega_f} q_t \partial_t((\delta_{jk} - a_{jk})\partial_j v_k)\right]_\tau^t$$
$$\lesssim_\lambda \left(\int_\tau^t \left(X^{3/2} + hX\right) + X^{3/2}(t) + hX(t) + X^{3/2}(\tau) + hX(\tau)\right) \lesssim hO(X),$$

where we used (3.15) and (3.28) in the last line.

**Step 3.** We show (3.21).

For $L_{\partial'\partial_t}$ and $L_{\partial_t}$ we use (3.15) and (3.26) to obtain

$$L_{\partial'\partial_t} = \lambda \int_{\Omega_f} \partial' q_t \partial' \operatorname{div} v = \lambda \int_{\Omega_f} \partial' q_t \partial'(\delta_{jk} \partial_j v_k) = \lambda \int_{\Omega_f} \partial' q_t \partial'((\delta_{jk} - a_{jk})\partial_j v_k)$$
$$\lesssim \lambda \|\partial' q_t\| \|I - a\|_2 \|v\|_2 \lesssim hX,$$

and similarly for $L_{\partial_t}$.

On the other hand, for the nonlinear part corresponding to $S \in \{\partial_t, \partial'\partial_t\}$, we have

$$N_{\partial'\partial_t} = \lambda \int_{\Omega_f} \partial' \partial_t((\delta_{jk} - a_{jl}a_{kl})\partial_k v_i)\partial_j \partial' v_i - \lambda \int_{\Omega_f} \partial' \partial_t((\delta_{ki} - a_{ki})q)\partial_k \partial' v_i$$
$$+ \int_{\Omega_f} (\delta_{jk} - a_{jl}a_{kl})\partial_k \partial_t \partial' v_i \partial' \partial_t v_i$$
$$\lesssim \|\partial_t((I - aa^T)\nabla v)\|_1 \|v\|_2 + \|\partial_t((I - a)q)\|_1 \|v\|_2 + \|I - aa^T\|_2 \|v_t\|_2^2 \lesssim hO(X),$$

due to (3.15), (3.26), and (3.27). In a similar way we obtain that $N_{\partial_t} \lesssim hO(X)$.

As for the cases $S \in \{\partial', \partial''\}$, we first note that

$$\|\phi\|_1 \leq \|S\eta(t) - S\eta(\tau)\|_1 + \|\widetilde{Sw}(\tau)\|_1 \lesssim \int_\tau^t \|v\|_3 + \|\overline{w}(\tau)\|_3 \leq O(X^{1/2}), \tag{3.29}$$

where we used (3.15) and $S\overline{w} = Sw$ on $\Gamma_c$ (recall (1.17)). Thus

$$N_S = \int_{\Omega_f} S((\delta_{jk} - a_{jl}a_{kl})\partial_k v_i)\partial_j \phi_i - \lambda \int_{\Omega_f} S((\delta_{ki} - a_{ki})q)\partial_k \phi_i + \int_{\Omega_f} (\delta_{jk} - a_{jl}a_{kl})\partial_k Sv_i \partial_j Sv_i$$
$$\lesssim \left( \|(I - aa^T)\nabla v\|_2 + \lambda \|(I - a)\nabla q\|_1 \right) \|\phi\|_1 + \lambda \int_{\Omega_f} Sa_{ki}q\partial_k \phi_i + \|I - aa^T\|_2 \|v\|_3^2$$
$$\lesssim ChO(X) + \lambda \int_{\Omega_f} Sa_{ki}q\partial_k \phi_i,$$

where we used (3.15), (3.26), (3.29) and we have singled out the most difficult term in which all derivatives (from $S$) fall onto $a$. In fact, since both $Sa$ and $q$ can accumulate from $t = 0$ (recall (2.24) and (3.26)), we can see that the last term is bounded by $h^2 O(X^{1/2})$, which is not sufficient for exponential decay (for which we need at least $hO(X)$, recall (3.23)). In order to obtain a stronger estimate we integrate by parts in $x_k$ and note that if the derivative falls on $q$ we can use (3.15) to estimate $\nabla q$, and if it falls on $a$, we can use the Piola identity (2.8). On the other hand, the boundary term can be estimated as in (2.20) to obtain additional $X^{1/2}$, since $S$ involves only horizontal derivatives. To be more precise, we obtain, in the case $S = \partial'$,

$$\int_{\Omega_f} \partial' a_{ki} q \partial_k \phi_i = -\int_{\Omega_f} \partial' a_{ki} \partial_k q \phi_i + \int_{\Gamma_c} \partial' a_{3i} q \partial' w$$
$$\lesssim \|a\|_2 \|\nabla q\|_1 \|\phi\|_1 + \|\partial' a_{3i}\|_{L^2(\Gamma_c)} \|q\|_{L^2(\Gamma_c)} \|\partial' w\|_{L^\infty} \lesssim hX,$$

where we used (2.8) and (3.5) in the first line, and (3.26), (3.15), (3.29), (2.24), (2.20) in the last inequality. The case $S = \partial''$ is similar, but we need to move one $\partial'$ away from $a$ in the boundary term. Namely,

$$\int_{\Omega_f} \partial'' a_{ki} q \partial_k \phi_i = -\int_{\Omega_f} \partial'' a_{ki} \partial_k q \phi_i + \int_{\Gamma_c} \partial' a_{3i} \partial'(q \partial' w)$$
$$\lesssim \|a\|_2 \|\nabla q\|_1 \|\phi\|_1 + \|\partial' a_{3i}\|_{L^2(\Gamma_c)} \left( \|\nabla q\|_{L^2(\Gamma_c)} \|\partial' w\|_{L^\infty} + \|q\|_{L^\infty} \|\partial'' w\|_{L^2(\Gamma_c)} \right)$$
$$\lesssim hX + X^{1/2}(\|\nabla q\|_1 \|\overline{w}\|_3 + \|q\|_2 \|\overline{w}\|_3) \lesssim hX,$$

as required, where we used (3.26), (3.15), and (2.20) in the second inequality and (3.15), (2.24) in the last.

For the commutator terms $C_S$, we recall (3.11) that

$$C_S = \int_{\Omega_f} (S\partial_j(a_{jl}a_{kl}\partial_k v_i) - \partial_j(a_{jl}a_{kl}\partial_k Sv_i))Sv_i$$
$$- \int_{\Omega_f} (S(a_{ki}\partial_k q) - a_{ki}\partial_k Sq)Sv_i + \int_{\Omega_f} Sq\,(a_{ki}S\partial_k v_i - S(a_{ki}\partial_k v_i)).$$

We observe that at each of the terms above at least one derivative from $S$ falls on $a$, and so we can use (3.15) and the fact that $\|a\|_2, \|a_t\|_2, \|a_{tt}\|_1 \lesssim h$ (a consequence of (3.26)) to obtain

$$\int_\tau^t C_S \lesssim hO(X)$$

if $S \in \{\partial_t, \partial', \partial'\partial_t, \partial_{tt}\}$. In the case $S = \partial''$ we obtain the same estimate, except in the cases when both derivatives fall on $a$ in the first integral. Indeed, this would result in three spatial derivatives falling on $a$, while we have no control over $\|a\|_3$. In this case we need to move one of $\partial'$ away from $a$. In other words, we obtain the term

$$\int_{\Omega_f} \partial_j(\partial''(a_{jl}a_{kl})\partial_k v_i)\partial'' v_i = \int_{\Omega_f} \partial''\partial_j(a_{jl}a_{kl})\partial_k v_i \partial'' v_i + \int_{\Omega_f} \partial''(a_{jl}a_{kl})\partial_k \partial_j v_i \partial'' v_i$$
$$= -\int_{\Omega_f} \partial'\partial_j(a_{jl}a_{kl})\partial'(\partial_k v_i \partial'' v_i) + \int_{\Omega_f} \partial''(a_{jl}a_{kl})\partial_k\partial_j v_i \partial'' v_i$$
$$\lesssim \|aa^T\|_2 \left(\|\nabla v\, \partial'' v\|_1 + \|D^2 v\|_{L^4}^2\right) \lesssim hX,$$

as required, where we used (3.26), (3.15), the inequality $\|\nabla v\|_{L^\infty} \lesssim \|v\|_3$ and the embedding $H^1 \subset L^4$.

## 3.3 An ODE-type lemma

In this section we prove an ODE-type result which is used to obtain the global bounds from (3.23). In what follows, the symbol $C$ denotes a constant that *does not change* its value from line to line.

**Lemma 4** (An ODE-type lemma). *Given $C \geq 1$, $\gamma \in (0,1]$, and $\lambda \in (0, 1/480C^2]$ there exists $\varepsilon > 0$ with the following property. Suppose that $f \colon [0,\infty) \to [0,\infty)$ is a continuous function satisfying*

$$f(t) + \lambda \int_\tau^t f \leq C(1 + \lambda^2(t-\tau))f(\tau) + C\lambda^2 \int_\tau^t f + h(t)O(f), \tag{3.30}$$

*for all times $t > 0$ such that*

$$h(t) := \sup_{(0,t)} f + \int_0^t f \leq \gamma, \tag{3.31}$$

*and all $\tau \in [0,t]$, where $O(f)$ denotes any term involving powers of $\lambda$, $C$ and $(t-\tau)$ and at least two factors of the form $f(\tau)$, $f(t)$ or $\int_\tau^t f$. Then the condition $f(0) \leq \varepsilon$ implies that*

$$f(t) \leq A\varepsilon e^{-t/a}, \tag{3.32}$$

*for all $t \geq 0$, where $a := 2C/\lambda$ and $A := 30C$.*

*Proof of Lemma 4.* First we note that letting $\varepsilon > 0$ be smaller than $\gamma/30C$ ensures that (3.31) holds for $t \in [0, T)$ for some $T > 0$, by continuity. Moreover, (3.31) holds at least as long as the claim (3.32) of the lemma holds. It is thus sufficient to verify (3.32), assuming that (3.31) holds for all times.

Suppose that the claim is false, let $t > 0$ be the first time such that
$$f(t) = A\varepsilon e^{-t/a},$$
and let $\tau \in [0, t)$ be the last time such that
$$f(\tau) = 2\varepsilon e^{-\tau/a};$$
note that $\tau > 0$, by continuity. As pointed out above, the inequalities (3.30) and (3.31) hold on $[0, t]$. We first claim that
$$t - \tau \leq \frac{4C}{\lambda}. \tag{3.33}$$
If the claim is false, then applying (3.30) on $[\tau, \tau + a]$, and taking into account only the second term on the left-hand side, gives
$$4C\varepsilon e^{-\tau/a} \left(1 - e^{-2}\right) = 2\varepsilon\lambda \int_\tau^{\tau+4C/\lambda} e^{-s/a}ds \leq \lambda \int_\tau^{\tau+4C/\lambda} f$$
$$\leq 2C\varepsilon(1 + 4C\lambda)e^{-\tau/a} + 4C^2\lambda A\varepsilon e^{-\tau/a} + \widetilde{C}\varepsilon^{3/2}e^{-\tau/a},$$
where the constant $\widetilde{C} > 0$ depends only on $C$, $\lambda$, and the form of the cubic term, namely the last term in (3.30). Noting that $1 - e^{-2} \geq 3/4$, and that $\lambda$ is sufficiently small so that $4C\lambda \leq 1/4$ and $120C^2\lambda \leq 1/4$, we may divide both sides of the above inequality by $\varepsilon e^{-\tau/a}$ to obtain
$$3C \leq \frac{5C}{2} + \frac{C}{4} + \widetilde{C}\varepsilon^{1/2},$$
which gives a contradiction if $\varepsilon = \varepsilon(\widetilde{C}) > 0$ is chosen sufficiently small. Thus (3.33) holds.

Applying (3.30) on $[\tau, t]$ and ignoring the second term on the left-hand side gives
$$30C\varepsilon e^{-t/a} = f(t) \leq 2\varepsilon C(1 + 4C\lambda)e^{-\tau/a} + 120C^3\lambda\varepsilon e^{-\tau/a} + \widetilde{C}\varepsilon^{3/2}e^{-\tau/a}.$$
Since $e^{-\tau/a} = e^{-t/a}e^{(t-\tau)/a} \leq e^{-t/a}e^2 \leq 10e^{-\lambda^{1/4}t/2}$, the smallness of $\lambda$ gives that $4C\lambda \leq 1/4$ and $120C^2\lambda \leq 1/4$ (as above). We can thus divide both sides by $\varepsilon e^{-t/a}$ to obtain
$$30C \leq 25C + \frac{5C}{2} + \widetilde{C}\varepsilon^{1/2},$$
which gives a contradiction for $\varepsilon = \varepsilon(\widetilde{C}) > 0$ sufficiently small. $\square$

## Acknowledgments

The authors are grateful to the referee for the careful reading of the manuscript and valuable comments and to Amjad Tuffaha for useful discussions. IK was supported in part by the NSF grants DMS-1907992 and DMS-2205493, while WSO was supported in part by the Simons Foundation.

I. Kukavica
Department of Mathematics, University of Southern California, Los Angeles, CA 90089
e-mail: kukavica@usc.edu

W. S. Ożański
Department of Mathematics, Florida State University, Tallahassee, FL 32304
and Institute of Mathematics, Polish Academy of Sciences, Warsaw 00-656, Poland
e-mail: wozanski@fsu.edu